\documentclass[11pt]{article}
\usepackage{amsmath}
\usepackage{amssymb}
\usepackage{amsthm}
\usepackage{graphicx}
\usepackage[round, authoryear]{natbib}
\usepackage{multirow}
\usepackage{tabularx}
\usepackage{longtable}
\usepackage[hidelinks]{hyperref}
\usepackage[left=1.25in,top=1.25in,right=1.25in,bottom=1.25in]{geometry}
\usepackage{fancyhdr}
 
\newtheorem*{theorem*}{Theorem}
\title{Motivated proofs: what they are, why they matter and how to write them}
\author{Rebecca Lea Morris}

\begin{document}

\maketitle
\thispagestyle{fancy}
\begin{abstract}
 Mathematicians judge proofs to possess, or lack, a variety of different qualities, including, for example, explanatory power, depth, purity, beauty and fit.  Philosophers of mathematical practice have begun to investigate the nature of such qualities.  However, mathematicians frequently draw attention to another desirable proof quality: being motivated.  Intuitively, motivated proofs contain no ``puzzling'' steps,  but they have received little further analysis.  In this paper, I begin a philosophical investigation into motivated proofs.  I suggest that a proof is motivated if and only if mathematicians can identify (i) the tasks each step is intended to perform; and (ii) where each step could have reasonably come from.  I argue that motivated proofs promote understanding, convey new mathematical resources and stimulate new discoveries.  They thus have significant epistemic benefits and directly contribute to the efficient dissemination and advancement of mathematical knowledge.  Given their benefits, I also discuss the more practical matter of how we can produce motivated proofs.  Finally I consider the relationship between motivated proofs and proofs which are explanatory, beautiful and fitting.
\end{abstract}

\tableofcontents

\section{Introduction}\label{Intro}
Mathematicians judge proofs to possess, or lack, a variety of different qualities.  For example, they may praise a proof for being explanatory, deep, pure, beautiful or fitting.  Philosophers of mathematical practice have attempted to clarify these virtues and explore their benefits.  With respect to mathematical explanation, \citet{steiner1978} and \citet{Kitcher1989}, for example, have proposed different theories.  On Steiner's approach, an explanatory proof is one which refers to a ``characterizing property'' and is, moreover, generalizable \citep[143--144]{steiner1978}.  However, according to Kitcher, an explanatory proof is one which is a member of the most unifying systematization of a set of mathematical beliefs \citep[\S 4]{Kitcher1989}.\footnote{Kitcher's main focus was scientific explanation, but he explicitly noted that his theory applied to mathematics as well \citep[437]{Kitcher1989}.}  The explanatory power of a proof may also connect directly to its depth.  For example, Lange has argued that ``depth is constituted (at least partly) by explanatory power'' \citep[footnote 1]{LangeDepth2015}.  On the topic of proof purity, \citet{Purity} have offered an analysis of what they call ``topical purity,'' and have further argued that proofs which are pure in this sense have particular epistemic benefits.  With respect to mathematical beauty, \citet{Rota} and \citet{Cellucci}, for example, have argued that beautiful proofs are ones which are ``enlightening'' or promote understanding.  Finally, \citet{Fit2016} have identified six senses in which a proof may be said to fit a theorem, and have further discussed how these senses of fit relate to other virtues, including beauty and explanation.

However, these are not the only qualities that mathematicians judge proofs for possessing or lacking.  In particular, mathematicians value proofs that are \emph{motivated}.  Intuitively, such proofs are perspicuous and do not contain any ``puzzling'' steps.  For example mathematicians have had the following to say about motivation (italics added): ``To someone unsteeped in Hermite's technique the \emph{motivation} for this proof must be unclear'' \citep[556]{IPI}; ``The most surprising aspect of this proof is the integral formulas, which have no apparent \emph{motivation}'' \citep[59]{CFE}; ``The proof of this generalization is simpler and the steps are better \emph{motivated} than the proof of the original result''\citep[4119]{HT}; ``A second purpose [of this paper] is \ldots to emphasize that this leads to much more clearly \emph{motivated} proofs''\citep[491]{rogers_lusins_1973}.\footnote{See \citep[\S 2.1]{Author2015} for more details and further examples.}

Moreover, mathematicians talk about motivated proofs in ways that suggest they have epistemic benefits.  P\'{o}lya's discussions, for example, indicate that unmotivated proofs have two important deficiencies: (i) they prevent the reader from fully understanding the argument; (ii) they fail to help the reader tackle her own mathematical problems \citep[\S 3]{WWM}.  Motivated proofs, on the other hand, are said to be free from such difficulties.  For example, when discussing two versions of the same proof, one with motivational material and the other without, P\'{o}lya wrote of the former ``Now we may understand how it was humanly possible to discover [a previously puzzling step] \ldots  The derivation \ldots becomes also more understandable'' \citep[690]{WWM}.  In other words, motivated proofs appear to help the reader better understand the argument and, by helping her to see how the steps could have been discovered, put her in a good position to adapt and reuse the ideas in her own mathematical work.

Although motivated proofs thus feature in mathematical practice and appear to have epistemic benefits, they have been subject to little analysis.  The goal of this paper, then, is to begin a philosophical investigation into motivated proofs.  I start in section (\ref{PM}) by discussing P\'{o}lya's work, which I refine and develop into an explication of motivated proofs in section (\ref{Def}).  In order to illustrate how this explication works in practice, I apply it to a proof of the Cauchy-Schwarz Inequality in section (\ref{CS}) and to a proof of the General Arithmetic-Geometric Mean Inequality in section (\ref{CS2}).  In section (\ref{value}) I analyze the epistemic benefits that motivated proofs provide, before addressing the more practical matter of how to produce motivated proofs in section (\ref{writing}).  In section (\ref{othervirtues}), I discuss the connection between motivated proofs and proofs which are explanatory, beautiful and fitting.  Finally, in section (\ref{conclude}) I call for an interdisciplinary investigation into motivational power in mathematics.

\section{P\'{o}lya on Motivated Proofs}\label{PM}
\subsection{``With, or Without, Motivation?''}\label{WWM}
Although there have been few attempts to clarify the notion of motivated proofs, P\'{o}lya's paper, ``With, or without, motivation?'', serves as a notable exception.  While P\'{o}lya did not present a precise definition of motivated proofs, he offered an example of a perfectly rigorous but unmotivated proof, discussed why it was deficient, and attempted to improve it with additional clarificatory material.\footnote{\citet{Sandborg1998, SandborgDiss} presented a nice discussion of P\'{o}lya's proof and subsequent analysis with respect to mathematical explanation.}  The unmotivated proof that P\'{o}lya analyzed was of a result called Carleman's inequality, which states the following (see e.g. \citep[684]{WWM}):

\begin{theorem*}[Carleman's Inequality]
Let $a_1, a_2, a_3, \ldots$ be a sequence of non-negative real numbers, not all equal to 0.  Then \[\sum_{n=1}^{\infty} (a_{1}a_{2}a_{3}\ldots a_{n})^{\frac{1}{n}} < e\sum_{n=1}^{\infty}a_{n}.\]
\end{theorem*}

To prove Carleman's inequality, we'll need to make use of the following three theorems, whose proofs I omit:\footnote{A proof of the General Arithmetic-Geometric Mean Inequality will be analyzed in section (\ref{CS2}).  For another proof and proofs of the other theorems, see e.g. \citep[20, 29, 30]{CSMasterclass}.}

\begin{theorem*}[Arithmetic-Geometric Mean Inequality (AGMI) ]
If $a_1$, $a_2$,\ldots, $a_n$ is a sequence of non-negative real numbers, then \[(a_{1}a_{2}\ldots a_{n})^{\frac{1}{n}}\leq \frac{a_1 + a_2 + \ldots + a_n}{n},\] with equality if and only if $a_1=a_2 = \ldots=a_n$.
\end{theorem*}

\begin{theorem*}[Sum of Telescoping Series (TS) ]
\[\sum_{n=k}^{\infty}\left(\frac{1}{n}-\frac{1}{(n+1)}\right)=\frac{1}{k}.\]
\end{theorem*}

\begin{theorem*}[Lower Bound for $e$ (LB) ]
For all $k = 1, 2, \ldots, n, \ldots$
\[\left(\frac{k+1}{k}\right)^{k}<e.\]
\end{theorem*}

We can now turn to P\'{o}lya's unmotivated proof of Carleman's inequality (see \citep[684-685]{WWM}):\footnote{I have edited P\'{o}lya's proof by including additional details to make it easier to follow.}

\begin{proof}
Define the sequence $c_1, c_2, c_3, \ldots $ by \[c_{1}c_{2}c_{3}\ldots c_{n}= (n+1)^n,\] for $n=1, 2, 3, \ldots$  Note that $c_{n}=\frac{(n+1)^{n}}{n^{n-1}}$.  We then have:

\begin{align*}
\sum_{n=1}^{\infty}\left(a_{1}a_{2}\ldots a_{n}\right)^{\frac{1}{n}} & = \sum_{n=1}^{\infty}\frac{\left(a_{1}c_{1}a_{2}c_{2}\ldots a_{n}c_{n}\right)^{\frac{1}{n}}}{n+1} \mbox{\phantom{.....}(def of $c_i$)} \\
& \leq \sum_{n=1}^{\infty}\frac{a_{1}c_{1} + a_{2}c_{2} + \ldots a_{n}c_{n}}{n(n+1)} \mbox{\phantom{.....}(AGMI)} \\
& = \sum_{k=1}^{\infty}a_{k}c_{k}\sum_{n\geq k}\frac{1}{n(n+1)} \\
& = \sum_{k=1}^{\infty}a_{k}c_{k}\sum_{n=k}^{\infty}\left(\frac{1}{n}-\frac{1}{n+1}\right) \\
& = \sum_{k=1}^{\infty}a_{k}\frac{(k+1)^{k}}{k^{k-1}}\cdot\frac{1}{k} \mbox{\phantom{.....}(def of $c_k$ \& TS)} \\
& = \sum_{k=1}^{\infty}a_{k}\left(\frac{k+1}{k}\right)^{k} \\
& < e\sum_{k=1}^{\infty}a_{k} \mbox{\phantom{.....}(LB)}
\end{align*}
\end{proof}

As previously mentioned, the above proof is perfectly rigorous, yet unsatisfactory because it fails to be motivated.  P\'{o}lya suggested that the motivational problem lies in the very first step: the introduction of the $c_i$ sequence.  As to why, exactly, it is problematic, he considered a variety of responses:

\begin{quote}
``It pops up from nowhere.  It looks so arbitrary.  It has no visible motive or purpose.''

\noindent``I hate to walk in the dark.  I hate to take a step, when I cannot see any reason why it should bring me nearer to the goal.''

\noindent``Perhaps the author knows the purpose of this step, but I do not and, therefore, I cannot follow him with confidence.''

\noindent[\ldots]

\noindent``Look here, I am not here just to admire you.  I wish to learn how to do problems by myself.  Yet I cannot see how it was humanly possible to hit upon your \ldots definition.  So what can I learn here?  How could I find such a \ldots definition by myself?'' \citep[685]{WWM}
\end{quote}

The problem with the $c_i$ sequence is thus \emph{not} that the reader cannot verify that it is used correctly, i.e. that it is used in accordance with the rules of logic.  Rather, the first three quotes suggest that the problem is that the imagined reader will find it difficult to identify its role in the argument.  The fourth quote suggests a related, though different problem: it is difficult for the intended readers to grasp the insight behind its introduction.\footnote{P\'{o}lya did not separate these two problems.  Instead, he suggested the problem was that the reader could not recognize the \emph{appropriateness} of the $c_i$ sequence, where ``[a] step of a mathematical argument is appropriate, if it is essentially connected with the purpose, if it brings us nearer to the goal'' \citep[685]{WWM}.}\textsuperscript{,}\footnote{There are strong similarities between P\'{o}lya's work and Mac Lane's.  For example, Mac Lane remarked ``[\ldots] there is some definite reason for the inclusion of each one of these steps in the proof; that is, each individual step is taken for some specific purpose'' \citep[125]{MacLane1935}.  He further explained ``To give a reason for a step means essentially to show why that step should be taken under the given conditions. In other words, given the form of the theorem to be proved, a knowledge of the steps already taken in its proof and a knowledge of the already established theorems which might be used in the next step, why should this particular step be taken?'' \citep[126]{MacLane1935}.}

Having diagnosed the problem, P\'{o}lya set out to fix it by providing a rational reconstruction of how the proof was discovered \citep[686-690]{WWM}.  This reconstruction helps readers to better grasp (i) the role that the $c_{i}$ sequence plays in establishing the theorem; and (ii) the insight that led to its introduction.\footnote{There may be a connection linking the distinction between the role of a step and the insight behind it and the distinction between ``how'' (or ``what'') and ``why.''  For example we might think of the role of a step as the answer to the how-question ``How does this step work?'' or the what-question ``What does this step do?''  Further, we might think of the insight behind the step as the answer to the why-question ``Why did the author take this step?''  However the connection may not be completely clear cut.  For instance, the insight behind a step could be considered to be the answer to the how-question ``How did the author come up with this step?'' or the what-question ``What led the author to introduce this step?'' rather than a why-question.   Similarly, although the why-question ``Why did the author take this step?'' can be naturally interpreted as inquiring after the insight behind a step, it could be interpreted as asking about the role that it plays instead.  That is to say, someone asking this question could be expecting an answer of the form ``The author took this step because it plays role $X$.''}   Below I will sketch part of P\'{o}lya's reconstruction to bring the role of the sequence and the insight that prompted its introduction into sharper focus.

First, it helps to know a general heuristic for working with inequalities, which I'll call Steele's heuristic: ``Wherever we hope to apply some underlying inequality to a new problem, the success or failure of the application will often depend on our ability to recast the problem so that the inequality is applied in one of those pleasing circumstances where the inequality is sharp, or nearly sharp'' \citep[26]{CSMasterclass}.  In other words, if we are going to apply an inequality to try to prove a new result, it's often a good idea to try to apply it when the inequality becomes (close to) an equality.

Now let's return to P\'{o}lya's proof.  He begins by noting that we are led to Carleman's inequality when we try to prove a related result, namely that $\sum_{n=1}^{\infty}a_{n} < \infty$ implies $\sum_{n=1}^{\infty}\left(a_{1}a_{2}\ldots a_{n}\right)^{\frac{1}{n}}<\infty$.  When attempting to prove this, a natural first step would be to apply the Arithmetic-Geometric Mean Inequality to the sequence $a_{1}, a_{2}, \ldots $.  This gives us: \citep[687]{WWM}:
\begin{align*}
  \sum_{n=1}^{\infty}\left(a_{1}a_{2}\ldots a_{n}\right)^{\frac{1}{n}} & \leq \sum_{n=1}^{\infty}\frac{a_{1} + a_{2} + \ldots + a_{n}}{n} \\
   & = \sum_{k=1}^{\infty}a_{k}\sum_{n=k}^{\infty}\frac{1}{n}
\end{align*}

As the last series is divergent, our proof attempt gets us nowhere.  The problem is that the $a_{i}$ sequence we applied the inequality to may have terms which are ``very unequal'' \citep[687]{WWM}.\footnote{The terms are unequal because $\sum_{n=1}^{\infty}a_{n}$ is assumed to be convergent.}  Why is this problematic?  Because the Arithmetic-Geometric Mean Inequality will be (nearly) sharp when and only when the terms of the sequence it is applied to are (nearly) equal.  In other words, applying the inequality directly to the $a_i$ sequence is not applying it where it is most effective, and so it is not surprising that the approach fails.  However, now that we know what the problem is, we can start to think of ways to solve it.  One straightforward way is to introduce a ``fudge factor'' \citep[28]{CSMasterclass} to make the terms of the sequence we apply the Arithmetic-Geometric Mean Inequality to more equal.  Thus we are led to the idea of introducing an auxiliary sequence, $c_i$.

The above sketch thus sheds some light on the insight that led to the introduction of the $c_i$ sequence, as well as the role that it plays in establishing the theorem.  In terms of insight, it reveals how Steele's heuristic prompted the introduction of an auxiliary sequence.  In terms of the role it plays, the sketch reveals that the sequence was chosen to make the terms of the $a_i$ sequence more equal.

However, there is more to be said about the insight behind the $c_i$ sequence and the role that it plays in the argument.  In particular, the sketch I presented focused only on the general strategy of introducing an auxiliary sequence---it did not say anything about the \emph{particular} choice of such a sequence.  P\'{o}lya's full rational reconstruction addresses this as well, but the considerations are more technical and so I omit them here.

\subsection{Discussion}\label{PDiss}
P\'{o}lya's analysis thus highlights that we want two things from motivated proofs: (i) to be able to identify the role of each step; (ii) to be able to identify the insight behind each step. In section (\ref{Def}) I will develop these desiderata into an explication of motivated proofs before illustrating how it is to be understood in practice by applying it to particular case studies in sections (\ref{CS}) and (\ref{CS2}).  To reflect the fact that different mathematical agents will differ in their ability to identify the role of a proof step or the insight behind it, the account that I offer will be a relative one.  More precisely, it will be relative to a given mathematical background, which I term the context.

After elaborating on the epistemic benefits of motivated proofs in section (\ref{value}), I will then consider ways in which we can write better motivated proofs in section (\ref{writing}).  In particular, while P\'{o}lya's rational reconstruction was an effective way to motivate the unsatisfactory proof of Carleman's inequality, I will show that motivational efficacy can often be improved by making more subtle changes to the proof itself.  As we shall see, issues relating to how a proof manages information are particularly important.

Before offering an explication of motivated proofs, however, I want to address potential objections to desideratum (ii) based on its connection to mathematical \emph{discovery}.  For example, it might be thought that the discovery process cannot be analyzed, and so we cannot say much about whether a mathematical agent reading a proof can identify the insight behind a given step.  However, as we have already seen from P\'{o}lya's rational reconstruction, it is often possible to analyze what lies behind these discoveries, and, as I will suggest in sections (\ref{Def})-(\ref{CS2}), it is also possible to assess whether an agent with a given mathematical background could reconstruct such reasoning.

Another potential objection to desideratum (ii) is that, being closely connected to discovery, it goes beyond the purview of philosophy even if it can be analyzed.  However, as Rav has emphasized, a core function of proofs is to convey new mathematical resources: ``Proofs are for the mathematician what experimental procedures are for the experimental scientist: in studying them one learns of new ideas, new concepts, new strategies---devices which can be assimilated for one's own research and be further developed'' \citep[20]{Rav1999}.  Part of doing this successfully involves not just conveying the resources themselves, but also information about when it is useful to try applying them.  Proofs that help readers to identify the insight behind each step provide more useful information about the circumstances under which their resources can be applied.  Thus desideratum (ii) connects closely to a core function of proofs and so should not be dismissed.

\section{An Explication of Motivated Proofs}\label{Def}
In this section I will sketch an explication of motivated proofs that incorporates the desiderata emphasized by P\'{o}lya.  I will build up to the notion of motivated proofs from the notion of motivated steps.  The account that I provide is thus reductionist.  However there may be other kinds of motivation that cannot be captured by a reductionist analysis and instead require a more holistic approach.  Such types of motivation may focus on, for example, a proof's interestingness or its connections to other results.\footnote{I am grateful to an anonymous reviewer for highlighting the possibility of other more holistic notions of motivation.}  These types of motivation are important and should also investigated, though they will not be my focus here.

Now for some preliminaries.  In what follows, when I refer to a \emph{proof}, I am referring to the proof as it is \emph{written}, for example, in a journal article.  I will assume that the reader of the proof is a \emph{finite} mathematical agent and, consequently, has only limited cognitive resources.  Moreover, while finite agents may differ in the amount of cognitive resources they have available, I will ignore these differences and focus instead on a typical finite agent, whose cognitive resources are neither impoverished nor exceptional.

Further, I will assume that the reader has access to a \emph{specific} mathematical ``toolkit,'' which I will call the \emph{context}.  The context contains mathematical resources, including definitions, theorems, proofs, techniques and heuristics.  By techniques, I mean precise approaches that can be applied to a variety of problems, for example the technique of completing the square in algebra, mathematical induction in number theory, or double counting in combinatorics.  By heuristics, I mean guidance about how to try solving problems, like Steele's heuristic.

Additionally, each resource comes with a \emph{range of application} within which it can be reliably applied.\footnote{I am intending ``reliably applied'' to be understood broadly, so as to include, for example, recognizing when a resource is, or can be, used.}\textsuperscript{,}\footnote{The context could, conceivably, contain more information.  It might include, for example, certain shared conventions among practitioners, or indications that certain resources are to be preferred over others.}  For example, a given resource may have a narrow range of application, so that it can only be reliably applied in very limited circumstances, whereas a resource with a wide range could be reliably applied even in new and unusual circumstances.  A medium range resource, then, would fall somewhere in between.

For instance, a technique concerning polynomials that has a narrow range of application might only be reliably applied in situations which are explicitly about polynomials.  However, a technique with a wide range of application could be reliably applied even if polynomials appear in a non-standard or unusual form.  Finally, if the technique has a medium range then it could be reliably applied even if it is not made explicit that polynomials are under discussion, but not if the polynomials appear in a non-standard or unusual form.  I will illustrate these considerations with respect to a particular case study in section (\ref{CS}).

Finally, in what follows I will only be considering proofs relative to contexts that allow typical finite agents with access to them to verify that the proof is correct.  This is because, if a typical finite agent with access only to a given context cannot verify that the proof is correct, she cannot recognize it as a proof and so it cannot be judged to be a motivated (or indeed unmotivated) proof relative to that context.

\subsection{Motivated Steps}
Let's say that an agent can identify the role a step plays in a proof if and only if she can identify the tasks that it is intended to perform.  Intended tasks can be purely logical, i.e. instantiating logical inference rules.  However, they often go beyond this and include things like introducing useful notation, concepts or definitions, simplifying expressions, instantiating mathematical resources or ``fine-tuning'' them to make the argument work in just the right way.  The task that the $c_i$ sequence was intended to perform, making the terms of the sequence the Arithmetic-Geometric Mean Inequality was applied to closer together, is an example of such ``fine-tuning.'' Now let's say that a step is \emph{role motivated} relative to a context $C$ if and only if a typical finite agent can identify the role that it plays, i.e. if and only if she can identify the tasks that it is intended to perform.\footnote{Role motivation can come in degrees.  For example, an agent may be able to identify, in a general way, a step's intended tasks but be unable to provide full details.  In such a case, the step would be partially role motivated relative to the context. }

Further, let's say that an agent can identify the insight behind a proof step if and only if she can identify where it could have reasonably come from.  Sometimes steps will come from the application of a mathematical technique.  For example, if we are using the technique of mathematical induction, we know that after we've proven the base case, the next step is to assume that the induction hypothesis holds for $n=k$.  In other words, that step comes directly from the application of the technique of induction.  Steps can also come from the application of a heuristic, like, for example, how the $c_i$ sequence in P\'{o}lya's proof came from the application of Steele's heuristic.  Sometimes steps come from chains of informal reasoning which can be complicated and rely on a variety of mathematical resources.  Now let's say that a step is \emph{insight motivated} relative to a context $C$ if and only if a typical finite agent with access to $C$ can identify the insight behind it, i.e. she can identify where it could have reasonably come from.\footnote{Insight motivation can come in degrees.  For example, an agent may be able to identify, in a general way, where a step could have reasonably come from, but be unable to provide full details.  In such a case, the step would be partially insight motivated relative to the context.}

Sometimes it is very easy for a typical finite agent with access only to a given context to identify a step's intended tasks or where it could have reasonably come from.  For example, a step's only intended task might be to instantiate part of a logical inference rule in a routine way.  Or it might come from the straightforward application of an available mathematical technique.  Other times, however, it requires considerably more effort to identify a step's intended task or where it came from.  For example, in P\'{o}lya's proof of Carleman's inequality, the intended task of the $c_i$ sequence and where it came from are only revealed when we think about what happens if we try to get by without introducing it and analyze what happens---the proof attempt fails because the terms of the $a_i$ sequence are too far apart and need to be made closer together so that the Arithmetic-Geometric Mean Inequality is applied when it is (nearly) sharp.

Next let's consider what factors affect whether a typical finite agent with access to a given context can identify the tasks a step was intended to perform or where it could have reasonably come from.  The resources that are involved in the step's intended tasks or where it reasonably came from need to be available in the agent's context and she must successfully connect them to the proof step in question for it to be role and insight motivated.  If the necessary resources are not available, then she cannot make the connection and so the step will fail to be role or insight motivated.  If they are available but only have a narrow range and the step's intended task or where it comes from involve the resources in a non-standard way, then the agent will be unable to make the connection and so the step will fail to be role or insight motivated.

Now suppose that all the resources a typical finite reader needs to identify the step's intended task and where it could have reasonably come from are contained within the context with wide enough range.  Nonetheless her finite nature will affect whether she can successfully identify the step's intended tasks or where it came from.  Making connections between the proof and the resources within her context requires cognitive effort, which she has in only limited supply.  If the cost of making the connection is too high, then our finite agent will not be able to make it and the step will fail to be role or insight motivated.  For example, a connection may be difficult to make, perhaps requiring the agent to make other connections first, and that is why its cost is too high.  Alternatively, the cost of the connection may be too high because the proof is poorly organized or includes a slew of irrelevant information, which makes the cost of making connections higher than it needs to be.  I will return to these factors when discussing how to write motivated proofs in section (\ref{writing}).

We can also consider how role and insight motivation are related.  In the discussion of the proof of the Generalized Arithmetic-Geometric Mean Inequality in section (\ref{CS2}) we will see that it is sometimes possible for a step to be role but not insight motivated.  Further, in the discussion of the Cauchy-Schwarz Inequality in section (\ref{CS}), we will see that it is sometimes possible for the insight behind a step to reveal its role.  This raises two follow-up questions: (i) is it sometimes possible for a step to be insight but not role motivated?; (ii) is it sometimes possible for the role of a step to reveal the insight behind it?\footnote{I am grateful to an anonymous reviewer for raising these questions.}  I do not have a definite answer to these questions, but in both cases I suspect the answer is ``yes.''

First consider question (i).  Although we will see that sometimes the insight behind a step involves its intended tasks and thereby reveals its role, it doesn't seem necessary for the insight behind a step to involve its intended tasks.  For instance, perhaps the insight behind a particular step focuses solely on what the proof looks like before its introduction, whereas its intended tasks only come into play in the later parts of the proof.  In such a case the insight behind the step will not, in and of itself, reveal its role.  Thus a reader with acess to a suitable context could identify the insight behind the step, i.e. identify where it could have reasonably come from, without identifying its role, i.e. its intended tasks.  In such a situation, the step would be insight but not role motivated.

Now consider question (ii).  When we're part way through proving a theorem and thinking about the next step to take, we might notice that if a certain task or tasks can be completed, then the proof will also be complete.  We can then use this information to (try to) find and introduce a step that achieves these tasks.  In such cases, the role of a step plays a central part in revealing where the step came from, i.e. the insight behind it.  Thus if a reader can identify the role of the step in such a case, so long as she has access to an appropriate context, she should also be able to identify the insight behind it.

Admittedly, some amount of vagueness remains in the notions of role and insight motivation.  This is to be expected as they must be open ended to capture a wide variety of tasks a step can be intended to perform and places a step could have reasonably come from.  However the case studies presented in sections (\ref{CS}) and (\ref{CS2}) will provide further illustrations of the terms and how they are to be applied in practice.

\subsection{Motivated Proofs}

We can use the notions of role and insight motivation, relative to a particular context, to sketch an explication of what it means for a proof to be motivated: a proof is motivated with respect to context $C$ if and only if each of its steps are both role and insight motivated relative to $C$.\footnote{Note that, just as role and insight motivation can come in degrees, so too can the motivational power of a proof.}  In other words, a proof is motivated relative to $C$ if and only if a typical finite reader with access only to $C$ can identify the tasks each step is intended to perform and where each step could have reasonably come from.

This explication thus makes the intuition that a motivated proof is one that doesn't contain any puzzling steps more precise.  Specifically, the notions of role and insight motivation rule out two different ways in a which a proof step could be puzzling.  For instance, a step that fails to be role motivated will cause puzzlement by leaving the reader in the dark as to what the step is being used to accomplish.  For example, a step that introduces a clever construction that is not role motivated will elicit questions like ``I see \emph{that} this works, but what, exactly, is it \emph{doing}?''  Similarly a step that fails to be insight motivated will cause puzzlement by leaving the reader at a loss as to how the proof author could have come up with it.  For example, a step that sets a variable to be a certain value that makes the proof work out perfectly but which fails to be insight motivated will be unsatisfying and lead to questions like ``What could have led the proof author to choose \emph{that} value?''  The conditions of role and insight motivation are thus necessary to prevent a step from being intuitively puzzling, as can be further appreciated from the discussion of the proof of Carleman's Inequality as well as the case studies in sections (\ref{CS}) and (\ref{CS2}).

I do not, however, have an argument to show that the conditions of role and insight motivation are sufficient to prevent a step from being intuitively puzzling.  While the puzzling steps in the proofs I examine both in this paper and elsewhere \citep{Author2015} have failed to be either role or insight motivated, it may nonetheless be possible for a step to be both role and insight motivated yet still intuitively puzzling.  If this turns out to be the case, the definition of motivated proofs I have given above will not fully capture the intuitive notion of motivation.  To fix it, further conditions would need to be added to rule out the additional ways in which a proof step can be puzzling.  The definition proposed above is thus best understood as a conjecture about the nature of motivated proofs.  While it's possible that the definition is incomplete, it still provides a useful starting point from which to characterize and analyze motivation in mathematics.

\section{Case Study: The Cauchy Schwarz Inequality}\label{CS}
To illustrate and further clarify the explication of motivated proofs just presented, I will apply it to a case study from linear algebra: the Cauchy Schwarz Inequality.
\subsection{Proof of The Cauchy Schwarz Inequality}
Before stating the Cauchy Schwarz Inequality, I will first briefly review real inner product spaces.  Let $V$ be a real vector space and $\langle , \rangle$ a function from $V\times V$ to $\mathbb{R}$ satisfying the following conditions\footnote{For more details on real inner product spaces, see e.g. \citep[7]{CSMasterclass}. Note that I have condensed the 5 axioms from Steele's presentation to 3.} for all $x, y, z \in V$:
\begin{enumerate}
\item $0\leq \langle x,x \rangle $ with equality if and only if $x=0$.
\item $\langle y,x \rangle =\langle x,y \rangle$.
\item for any $\alpha, \beta \in \mathbb{R}$ $\langle \alpha x + \beta y, z \rangle = \alpha \langle x,z \rangle + \beta \langle y, z \rangle.$
\end{enumerate}
Then $(V, \langle , \rangle)$ is called a \emph{real inner product space.}  Further, given a real inner product space, we can define a new function, $|| \ . \ || : V \rightarrow \mathbb{R}$, called a \emph{norm}, as follows: $||x|| = \langle x,x \rangle^{\frac{1}{2}}$.

We can now state the Cauchy Schwarz inequality for a real inner product space:
\begin{theorem*}
Let $(V, \langle , \rangle)$ be a real inner product space.  Then for all $x, y \in V$
\[|\langle x,y \rangle| \leq ||x||\cdot||y||\] with equality if and only if $x$ and $y$ are linearly dependent.
\end{theorem*}

Consider the proof, below, which loosely follows the presentation by \citet[3]{deutsch}, though includes more details.
\begin{proof}
First suppose that $x$ and $y$ are linearly dependent.  Then $x = \alpha y$ for some $\alpha \in \mathbb{R}$.  Thus we have

\[|\langle x,y \rangle| = |\langle \alpha y, y \rangle| = |\alpha \langle y,y \rangle| = |\alpha|||y||^{2} = ||\alpha y||\cdot ||y|| = ||x||\cdot ||y||.\]

Now suppose that $x$ and $y$ are linearly independent.  Then for no $\alpha \in \mathbb{R}$ do we have $y-\alpha x = 0$.  Thus for all $\alpha$
 \begin{align*}
0< \langle y-\alpha x, y-\alpha x \rangle & = \langle y, y-\alpha x \rangle - \alpha \langle x, y-\alpha x \rangle \\
& = \langle y,y \rangle -\alpha \langle x, y \rangle -\alpha \langle y, x \rangle + \alpha^{2}\langle x,x \rangle \\
& = \langle y,y \rangle -2\alpha \langle x, y \rangle + \alpha^{2}\langle x,x \rangle .
\end{align*}

Now let $\alpha : = \frac{<x,y>}{<x,x>}$.  Substituting this in the above yields
\begin{align*}
0 < \ \langle y-\frac{\langle x,y \rangle}{\langle x,x \rangle}x, y -\frac{\langle x,y\rangle }{\langle x,x \rangle} x \rangle & = \langle y,y \rangle -2\frac{\langle x,y \rangle}{\langle x,x \rangle} \langle x, y \rangle + \frac{\langle x,y\rangle^{2}}{\langle x,x \rangle^{2}}\langle x,x\rangle \\
& = \langle y,y \rangle - \frac{\langle x,y\rangle ^{2}}{\langle x,x\rangle} .
\end{align*}
Rearranging thus yields
\[\frac{\langle x,y\rangle ^{2}}{\langle x,x\rangle} < \ \langle y,y \rangle.\]
Hence
\[\langle x,y \rangle ^{2} < \ \langle x,x\rangle \langle y,y\rangle.\]
Finally, taking square roots of both sides,
\[|\langle x,y \rangle | < ||x|| \cdot ||y||.\]
\end{proof}
\subsection{Motivational Efficacy}
I will focus on one step in the proof of the Cauchy Schwarz Inequality: the step in which $\alpha$ is set to be $\frac{<x,y>}{<x,x>}$.  As I will describe below, this step performs a specific task within the proof and comes from an application of Steele's heuristic.  Nonetheless, the step fails to be role and insight motivated, relative to a particular context that I will describe later on.  This is because, while a typical agent with access only to that context can verify the correctness of the proof, she cannot identify the step's intended task or where it could have reasonably come from.

First, to facilitate discussion, I am extracting and numbering the steps surrounding the choice of $\alpha$ below:

\begin{itemize}
\item[] For all $\alpha \in \mathbb{R}$, $0 < \langle y,y \rangle -2\alpha \langle x, y \rangle + \alpha^{2}\langle x,x \rangle$ \hspace*{\fill} (1)
\item[] Let $\alpha : = \frac{<x,y>}{<x,x>}$ \hspace*{\fill} (2)
\item[] Then $0 < \langle y,y \rangle - \frac{\langle x,y\rangle ^{2}}{\langle x,x\rangle}$ \hspace*{\fill} (3)
\end{itemize}

To start, let's focus on the intended task of step (2): to introduce a value of $\alpha$ that minimizes $\langle y,y \rangle -2\alpha \langle x, y \rangle + \alpha^{2}\langle x,x \rangle$.  That this value of $\alpha$ successfully performs this task can be checked using techniques from calculus or by completing the square, since $\langle y,y \rangle -2\alpha \langle x, y \rangle + \alpha^{2}\langle x,x \rangle$ is a quadratic in $\alpha$.  For example, here is Steele's demonstration via completing the square \citep[57]{CSMasterclass}:

\begin{align*}
\langle y-\alpha x,y-\alpha x\rangle & = \langle y,y \rangle -2\alpha \langle x, y \rangle + \alpha^{2}\langle x,x \rangle \\
& = \langle x,x \rangle \left(\frac{\langle y, y \rangle}{\langle x, x \rangle} -2\alpha\frac{\langle x,y \rangle}{\langle x,x \rangle} +\alpha^{2} \right) \\
& =  \langle x,x \rangle \left \{ \left(\alpha - \frac{\langle x,y \rangle}{\langle x, x \rangle}\right)^{2} + \frac{\langle y,y \rangle}{\langle x,x \rangle} - \frac{\langle x,y \rangle^{2}}{\langle x,x \rangle^{2}}\right \}
\end{align*}

In the above, $\langle x,x \rangle$ is non-negative and $\frac{\langle y,y \rangle}{\langle x,x \rangle} - \frac{\langle x,y \rangle^{2}}{\langle x,x \rangle^{2}}$ is a constant.  Thus to minimize $\langle y-\alpha x,y-\alpha x\rangle$, we need to minimize $\left(\alpha - \langle x,y \rangle / \langle x, x \rangle\right)^{2}$.  But this is always non-negative, taking its minimum value of 0 when $\alpha = \frac{\langle x,y \rangle}{\langle x,x \rangle}$.

Now let's consider where step (2) came from.  First recall that we have an inequality in step (1).  Then remember Steele's heuristic: it's generally a good idea to apply inequalities where they are (nearly) sharp.  This supports introducing a value of $\alpha$ that minimizes the right hand side of the inequality, as this will make the inequality as sharp as possible.

However, while typical agents with access only to the context $C$, described below, will be able to check the correctness of the proof, they will fail to identify step (2)'s intended task (introducing a value that minimizes a quadratic) and where it could have reasonably come from (Steele's heuristic).  Thus, relative to this context, step (2) fails to be both role and insight motivated and so the proof itself fails to be motivated.  Let context $C$ be described in broad terms as follows:

\begin{itemize}
\item[]\textbf{Medium range resources}: \emph{Real inner product spaces} (definition of a real inner product space, definition of norm, definition of vector space, definition of linear (in)dependence, etc); \emph{Algebra} (theorem stating that the min or max of a quadratic $ax^2 + bx + c$ occurs at $-b/2a$, technique of completing the square, techniques for manipulating algebraic expressions); \emph{Logic} (standard inference rules)
\item[]\textbf{Missing resources}: \emph{Steele's heuristic} (try to apply inequalities where they are sharp)\footnote{Note that this context is not ``artificial'' for lacking Steele's heuristic.  Steele's heuristic is somewhat technical and agents do not need to have it available to them in order to understand linear algebra.  So, while it may be included in sophisticated contexts, omitting in from an intermediate context like $C$ does not make $C$ unrealistic.}
\end{itemize}

Let's first consider whether the step in question is role motivated relative to $C$.  There are two ways a reader could identify the step's intended task.  The first is by connecting the expression $\langle y, y \rangle - 2\alpha\langle x, y\rangle + \alpha^2 \langle x,x \rangle$ from step (1) and the value $\frac{\langle x, y\rangle}{\langle x,x\rangle}$ from step (2) to the algebraic resources in her context.  For example, the reader could recognize that this pair have the general form $ax^2 + bx + c$ and $-b/2a$ and then make a connection to the theorem in her context that states that the minimum or maximum value of a quadratic $ax^2 + bx +c$ occurs at $-b/2a$.  This would then allow her to infer that the value of $\alpha$ is chosen to minimize $\langle y, y \rangle - 2\alpha\langle x, y\rangle + \alpha^2 \langle x,x \rangle$ and she would have successfully identified the step's intended task.

However, a typical finite agent with access to context $C$ cannot connect $\langle y, y \rangle - 2\alpha\langle x, y\rangle + \alpha^2 \langle x,x \rangle$ and the value $\frac{\langle x, y\rangle}{\langle x,x\rangle}$ to the quadratic resources within her context.  This is because these resources only come with a medium range, meaning they cannot be reliably applied (or recognized as applying) when quadratics appear in a non-standard or unusual form.  Yet the quadratic in the proof \emph{does} appear in a non-standard form.  More precisely, three features of the expression $\langle y, y \rangle - 2\alpha\langle x, y\rangle + \alpha^2 \langle x,x \rangle$ obscure its quadratic nature.  First, the expression uses a non-standard choice of variable, $\alpha$.  Second, the coefficients have a complex form, being expressed in terms of inner products, and are written to the right of the variable, not to the left as is standard.  Third, while quadratics are standardly written starting with the second degree term then the first degree term and finally the constant term, this one is written in reverse.  Consequently a typical finite agent cannot identify the step's intended task in this way.

The second way that a reader could come to identify the intended task of step (2) is by first identifying where it could have reasonably come from.  For example, if the reader recognizes that it is a good idea to exploit the inequality $0<\langle y,y\rangle -2\alpha\langle x,y\rangle + \alpha^2\langle x,x\rangle$ where it is as sharp as possible, she can then infer that this means minimizing the right hand side.  Next, she could use the algebraic resources in her context (either the technique of completing the square or the theorem concerning min/max values) to confirm that the choice of $\alpha$ is the minimizing value and thus identify the intended task of step (2).

However, a typical finite agent with access to context $C$ cannot identify step (2)'s intended task in this way.  First, Steele's heuristic, which states that it is useful to apply inequalities where they are as sharp as possible, is not available.  Thus the reader cannot complete the first part of the reasoning, i.e. cannot recognize that the proof should exploit the inequality  $0<\langle y,y\rangle -2\alpha\langle x,y\rangle + \alpha^2\langle x,x\rangle$ where it is as sharp as possible.  Second, the quadratic $\langle y,y\rangle -2\alpha\langle x,y\rangle + \alpha^2\langle x,x\rangle$ still appears in non-standard form, so the reader cannot reliably apply the quadratic resources from her context.  Thus the reader cannot complete the second part of the reasoning, either.

The considerations in the previous paragraph also show that step (2) fails to be insight motivated relative to $C$.  This is because, in order to identify where the step could have reasonably come from the reader needs to apply Steele's heuristic, but this resource is not available within her context.  Consequently, step (2) is neither role nor insight motivated relative to $C$.  Thus the proof fails to be motivated relative to $C$.

However, we can use our analysis of why the proof fails to be motivated relative to $C$ to improve it.  There are two issues that we need to address: (i) the contextual resources about quadratics only have a medium range of application; (ii) the contextual resources do not include Steele's heuristic.

Consider first the issue of the quadratic resources having only a medium range of application.  If the proof explicitly points out that $\langle y,y\rangle -2\alpha\langle x,y\rangle + \alpha^2\langle x,x\rangle$ is a quadratic in $\alpha$, its nature will no longer be obscured.  Thus such a comment will bring the quadratic within the range of application of the contextual resources, and so put the reader in a better position to successfully make connections between her context and the proof.  This will then allow her to identify the intended task of step (2).

Similarly, to compensate for the fact that Steele's heuristic is missing, the proof could include a small remark to explain that applying the inequality where it is as sharp as possible allows the most information to be extracted from it. The reader could then identify where this particular value of $\alpha$ could have come from.  Moreover, it would provide her with the heuristic which she could then reuse in future.

More concretely, we could replace the part of the original proof in which the particular value of $\alpha$ is introduced with the following:

\begin{quote}
To get the most information out of this inequality, we should apply it where it is as sharp as possible.  Thus we should apply it when the quadratic $\langle y, y\rangle - 2\alpha \langle x,y\rangle + \alpha^2 \langle x,x\rangle$ takes its minimum value.  So let $\alpha = \frac{\langle x, y\rangle}{\langle x, x\rangle}$ and substitute this back into the inequality \ldots
\end{quote}

This relatively minor tweak should be sufficient to significantly improve the proof's motivational power relative to context $C$.

\section{Case Study: The General Arithmetic-Geometric Mean Inequality}\label{CS2}
As another illustration of the explication of motivated proofs presented in section (\ref{Def}), I will now apply it to a proof of the General Arithmetic-Geometric Mean Inequality.

\subsection{Proof of the General Arithmetic-Geometric Mean Inequality}
The General Arithmetic-Geometric Mean Inequality can be stated as follows (see e.g. \citep[23]{CSMasterclass}):
\begin{theorem*}
Let $p_1, p_2, \ldots, p_n$ be non-negative real numbers such that $\sum_{i=1}^{n}p_{i}=1$ and $a_1, a_2, \ldots, a_n$ be non-negative real numbers.  Then we have $a_{1}^{p_{1}}a_{2}^{p_{2}}\ldots a_{n}^{p_{n}} \leq p_{1}a_{1}+p_{2}a_{2}+\ldots +p_{n}a_{n}$.
\end{theorem*}

In fact, when proving this we can assume all of the $p_k$ and $a_k$ are positive.  If $p_i=0$, for example, we can ignore the terms corresponding to this on either side of the inequality as they do not contribute anything.  If $p_i\neq 0$ and $a_i = 0$ then $a_{1}^{p_{1}}a_{2}^{p_{2}}\ldots a_{n}^{p_{n}} =0$ and the inequality holds.

A proof of this inequality, which P\'{o}lya discovered while dreaming, makes use of the following bound for the exponential function: $1+x \leq e^{x}$ for all $x\in\mathbb{R}$ with equality if and only if $x=0$.\footnote{For a variety of different proofs of this bound, see the answers to the Math StackExchange post by \citet{InequalityProofs}.}  By a change of variables from $x$ to $x-1$ this becomes:

\begin{theorem*}[Bound on $e^{x-1}$]
For all $x\in\mathbb{R}$, $x \leq e^{x-1}$ with equality if and only if $x=1$.
\end{theorem*}

With this theorem in hand, we can now prove the General Arithmetic-Geometric Mean Inequality as follows (see e.g. \citep[23--25]{CSMasterclass}, \citep{AoPS}).

\begin{proof}
Let $G=a_{1}^{p_{1}}a_{2}^{p_{2}}\ldots a_{n}^{p_{n}}$ and $A= p_{1}a_{1}+p_{2}a_{2}+\ldots +p_{n}a_{n}$.  Now define $\alpha_{k}=\frac{a_{k}}{A}$ for $1\leq k \leq n$.  We apply the bound on $e^{x-1}$ when $x=\alpha_{k}$ to get $\alpha_{k} \leq e^{\alpha_{k}-1}$, i.e. \[\frac{a_{k}}{A}\leq e^{\frac{a_{k}}{A}-1}.\]
Raising to the power $p_{k}$ and multiplying we get
\[\left(\frac{a_{1}}{A}\right)^{p_{1}}\left(\frac{a_{2}}{A}\right)^{p_{2}} \ldots \left(\frac{a_{n}}{A}\right)^{p_{n}} \leq \exp{\left(\left\{\sum_{k=1}^{n}p_{k}\frac{a_{k}}{A}\right\}-1\right)}=1.\]
Rearranging we obtain
\[a_{1}^{p_{1}}a_{2}^{p_{2}}\ldots a_{n}^{p_{n}} \leq A^{p_{1} +p_{2} + \ldots + p_{n}},\] i.e.
\[G \leq A^{p_{1} +p_{2} + \ldots + p_{n}}.\]
As $p_{1}+p_{2}+\ldots + p_{n} =1$ we thus have $G\leq A$

\end{proof}

\subsection{Motivational Efficacy}
As in the previous case study, I will focus only on one step in the proof of the General Arithmetic-Geometric Mean Inequality: the step in which $\alpha_{k}$ is defined to be $\frac{a_{k}}{A}$.  As I will describe below, this step performs a specific task within the proof and comes from recognizing that the general theorem can be reduced to a special case.  Relative to a particular context that I describe later on this step will be role, but not insight, motivated.  This is because, while a typical finite agent with access only to that context can verify the correctness of the proof and identify the step's intended task, she cannot identify where it could have reasonably come from.

First, to facilitate discussion, I am extracting and numbering the steps surrounding the definition of $\alpha_{k}$ below:

\begin{itemize}
\item[] Let $G=a_{1}^{p_{1}}a_{2}^{p_{2}}\ldots a_{n}^{p_{n}}$ and $A= p_{1}a_{1}+p_{2}a_{2}+\ldots +p_{n}a_{n}$. \hspace*{\fill} (1)
\item[] Define $\alpha_{k}=\frac{a_{k}}{A}$ for $1\leq k \leq n.$ \hspace*{\fill} (2)
\item[] Apply the bound on $e^{x-1}$ when $x=\alpha_{k}$ to get $\alpha_{k} \leq e^{\alpha_{k}-1}.$ \hspace*{\fill} (3)
\end{itemize}

The intended task of step (2) is relatively straightforward: to rescale the variables $a_{k}$ so that $\sum_{k=1}^{n}p_{k}\alpha_{k}=1$.  That the definition of $\alpha_{k}$ achieves this is clear by inspection.  The insight behind step (2), however, is more complicated.  Below I sketch Steele's reconstruction of how the proof could have been discovered, as this nicely illustrates where step (2) comes from \citep[23--25]{CSMasterclass}.

Steele starts by applying the bound for $e^{x-1}$ when $x=a_{k}$ to obtain $a_{k}\leq e^{a_{k}-1}$.  Raising each side to the power $p_{k}$ and then multiplying the resulting inequalities he obtains  \[G=a_{1}^{p_{1}}a_{2}^{p_{2}}\ldots a_{n}^{p_{n}} \leq \exp{\left(\left\{\sum_{k=1}^{n}p_{k}a_{k}\right\}-1\right)} = e^{A-1}.\]

By applying the bound for $e^{x-1}$ when $x=A$ he also finds \[A \leq e^{A-1}.\]

In other words, he has a double bound, expressed by \[\max{\{A, G\}} \leq e^{A-1}.\]

If $A=e^{A-1}$ then the General Arithmetic-Geometric Mean Inequality will be established.  $A=e^{A-1}$ when and only when $A=1$, so the inequality is in fact proven for the special case when $A=1$.  Steele now tries to reduce the general case to the special one.\footnote{``Try finding a special case and reduce the general case to the special one'' is a general problem solving heuristic.  \citet{ReduceToSimple} briefly discusses this heuristic in connection with two examples.}  This can be done by rescaling the $a_{k}$ so that their weighted sum is 1.  Hence the introduction of the variables $\alpha_{k}=\frac{a_{k}}{A}$ since then $\sum_{k=1}^{n}p_{k}\alpha_{k}=1$.

In summary, step (2) comes from recognizing that we can prove the theorem for a special case and that we can reduce the general case to the special one by rescaling variables.

Having seen both the intended task of step (2) (to rescale the $a_k$) and where it could have reasonably come from (reducing the general to a special case), let us consider whether the proof is motivated relative to the intermediate context $C$ described below.

\begin{itemize}
\item[]\textbf{Medium range resources}: \emph{Inequalities} (techniques for working with inequalities, e.g. rescaling variables, raising inequalities to powers, multiplying inequalities; theorem for bound on $e^{x-1}$); \emph{Algebra} (techniques for operating with sums and products); \emph{Problem solving} (heuristics such as ``Try finding a special case and reduce the general case to the special one''); \emph{Logic} (standard inference rules)
\end{itemize}

First note that a typical finite agent with access to context $C$ will have no trouble checking the correctness of the proof of the General Arithmetic-Geometric Mean Inequality. Moreover, such an agent can identify the intended task of step (2): to introduce rescaled variables whose weighted sum is 1.  This is because the technique of rescaling is available within the context with medium range, and there is nothing non-standard or unusual that would place step (2) outside of this range.  This means that step (2) is role motivated relative to $C$.

Step (2) fails to be insight motivated relative to $C$, however.  Recall that we saw, above, that step (2) comes from reducing the general case to a special one.  However, to reconstruct this, the agent must first realize that the theorem can be proven in a special case and then figure out how to reduce the general case to the special one.  Recognizing that the theorem can be proven in a special case itself involves a number of steps: Show $G\leq e^{A-1}$; Infer $A\leq e^{A-1}$; Infer $\max{\{A, G\}}\leq e^{A-1}$; Conclude $G\leq A$ if $A=1$.  In other words, the reasoning involved in identifying where step (2) comes from has a number of different pieces and the agent has to complete each of them and put them together successfully.  Further, the above proof of the inequality does not provide the reader with any prompts to help her complete these parts.  Rather, she has to complete them on her own, without guidance.  As she is only a typical finite agent, however, the more parts she must complete on her own like this, the less likely it is that she can complete them all successfully.  In other words, as the reasoning needed to identify where the step comes from involves numerous parts, it makes it less likely that the step will be insight motivated.

Moreover, although the heuristic ``Try finding a special case and reduce the general case to the special one'' is available in $C$, a typical finite agent won't be able to apply it in this case.  This is because once the agent has obtained the double bound $\max{\{A, G\}} \leq e^{A-1}$ she may feel as though she is at a dead end and not realize she has established the theorem for a special case.  Indeed Steele refers to the task of obtaining an inequality between $G$ and $A$ from a bound on their maximum as ``a modest paradox'' \citep[24]{CSMasterclass} and remarks that a reader ``might be discouraged'' \citep[25]{CSMasterclass} at this point.  If the reader does not recognize she has established the theorem for a special case, then she cannot try to reduce the general case to it.  So even if she can reconstruct the reasoning to this point, she will not be able to complete it and so step (2) will fail to be insight motivated relative to $C$.

As step (2) fails to be insight motivated, the proof fails to be motivated relative to $C$.  However, again our analysis can help us determine how to better motivate the proof.  There were two main problems for an agent trying to determine where step (2) came from: (i) she needed to reconstruct multiple pieces of reasoning entirely on her own, with no guidance from the proof; (ii) part of the reasoning looked like a dead-end.  These issues can be remedied by preceding the proof with a sketch of the reconstruction of its discovery.  In fact, this is exactly how \citet[23--25]{CSMasterclass} proceeds, so his version of the proof is motivated relative to $C$.

\section{Epistemic Benefits of Motivated Proofs}\label{value}
A proof that is motivated relative to a given context provides its readers (i.e. the practitioners of that context) with epistemic benefits.  In particular, taking `understanding' in a pre-theoretical sense, motivated proofs promote understanding, convey new mathematical resources and stimulate new discoveries.

Mathematicians often point to understanding as the main goal of mathematics.  For example, Thurston remarked ``The measure of our success is whether what we do enables people to understand and think more clearly and effectively about mathematics'' \citep[163]{Thurston}.  A proof that is motivated relative to a given context promotes understanding among practitioners of that context by providing them with more useful information than an unmotivated proof.  More precisely, a motivated proof helps practitioners to identify (i) the intended tasks of each step in the proof; and (ii) where each step could have reasonably come from.  This additional information enables readers to better understand how the different proof steps work together to establish the result.  Further, it promotes understanding of the theorem itself.  In particular, if a reader can identify the intended task of each step and where they reasonably came from, she should be able to grasp why certain conditions are included in the statement of the theorem, or whether they can be loosened, for example.\footnote{It is also possible for a reader to come to gain understanding of a theorem independently.  In particular, I do not mean to imply that a reader will only fully understand a theorem after reading a motivated proof.  A reader may gain understanding of a theorem by, for example, considering specific instances or special cases of the theorem or via visualization.}

Rav has drawn attention, as we saw in section (\ref{PDiss}), to the fact that a central function of proofs is to convey new mathematical resources.  Part of successfully conveying these resources means showing how they can be used and communicating when it is useful to apply them.  Proofs that are motivated relative to a given context will be more successful at this than unmotivated ones.\footnote{Motivated proofs thus have the instrumental benefit of promoting reuse of their resources.  See \citep{Author2019} for details about reuse in mathematics. }  In particular, as practitioners can identify the intended task of each step in a motivated proof, they will better learn how the resources within the proof can be used.  By identifying where each step could have reasonably come from, the practitioners will also be learning when it is useful to try to apply these resources in future.

Moreover, motivated proofs can sometimes directly stimulate discovery of new results.  For example, identifying a step's intended task or where it reasonably came from can suggest new proofs of new results.\footnote{This means that sometimes \emph{un-}motivated proofs can prompt discovery in this way because some (but not all) steps in an unmotivated proof may be role and insight motivated relative to the given context.  However, as \emph{all} steps in a motivated proof are role and insight motivated, they will prompt such discovery more reliably than unmotivated ones.}  For a concrete example of this, consider the following toy example: every odd integer can be represented as the difference of two squares.  The proof of this is very simple.  If $n$ is an odd integer, then we can write $n=2k+1$ for some integer $k$.  Then we just have to notice that $(k+1)^2 - k^2 = 2k+1$ and the proof is complete.

The only non-trivial part of this proof is the insight behind the introduction of the witnesses to the desired representation of $n$.  However, here's a rational reconstruction of where these witnesses came from (see e.g. \citep{Clark2Squares}): To prove the theorem, we want to find integers $x$ and $y$ such that  $2k+1 = x^2 - y^2.$  We can rewrite this as  $2k + 1 = (x-y)(x+y).$  Now, to try to make the left hand side of the above equation look more like the right hand side, we can again rewrite it as $1\cdot(2k+1) = (x-y)(x+y).$  At this point, \emph{we make a guess} to see if we can get further and try setting $x-y=1$ and $x+y=2k+1$.  Then we have a pair of simultaneous equations in two unknowns.  Moreover, we find that they are solvable, yielding $x=k+1, y=k$.

Notice that part of this informal reasoning involves a guess (which is italicised).  This may prompt us to ask: are there circumstances under which it \emph{isn't} a guess?  And the answer is: yes!  If $p$ is an odd \emph{prime} then its only divisors are 1 and itself.  Thus $(x-y)$ must be 1 and $(x+y)$ must be $2k+1$.  Recognizing this, we may then be led to a proof of a new result: every odd \emph{prime} can be represented \emph{uniquely} as a difference of two squares.  Notice that this does not hold for odd numbers generally since, for example, $15 = 8^2 - 7^2$ and $15 = 4^2 - 1^2$.

Motivated proofs thus provide significant benefits to the practitioners of a given context.  By promoting understanding, they help mathematicians to achieve one of their fundamental goals.  By conveying how and when to apply mathematical resources, they better serve one of the core functions of proof.  By directly stimulating discovery, they contribute to the advancement of mathematical knowledge.

\section{Writing Motivated Proofs}\label{writing}
As motivated proofs thus have significant benefits, we should consider the more practical matter of how we can produce them.\footnote{I am not claiming that it is always possible to write motivated proofs with respect to a given context.  Some proofs may only appear motivated to agents with access to suitably sophisticated contexts and perhaps there are proofs that are not motivated relative to any known context.}  To do this, we need to consider both the context that a proof will be assessed against, and the fact that readers of the proof are finite agents with limited cognitive resources.

In terms of the context, we need to consider whether the resources needed to identify each step's intended task and where each step could have reasonably come from are (i) available in the context; (ii) associated with a wide enough range of application.  If required resources are not present in the context, then, to produce a motivated proof, the writer will need to include additional content to make up for the missing resources.  If the required resources are present, but do not have a wide enough range of application, then the writer will need to help the reader connect the relevant resources to the proof.  The suggested improvements to the proofs of the Cauchy Schwarz Inequality discussed in section (\ref{CS}) are of this form, for example.

However, to ensure that we produce motivated proofs we also need to take into account the fact that the reader is a \emph{finite} mathematical agent, with limited cognitive resources.  Thus the proof writer should aim to make it as easy as possible for her to make connections between the resources in her context and the proof, even if the required resources are available and associated with a wide enough range of application.  This is because, the more effort it takes for a reader to make those connections, the more difficult it is to make them and the less likely it is that a typical reader can identify each step's intended task and where it reasonably came from.

If, for example, a lengthy chain of reasoning is needed to identify a step's intended tasks or where it reasonably came from, the reader must make all of the connections in the chain to do so, which will be cognitively expensive.  To address this, the proof writer could include prompts to make some of the connections for the reader, so that she does not need to make them all by herself.  Similarly, if a single connection is particularly difficult, i.e. cognitively expensive, to make, the proof writer should add a prompt to reduce the cost of making the connection.  The suggested improvements to the proof of the General Arithmetic-Geometric Mean Inequality discussed in section (\ref{CS2}) are of this form, for example.

Furthermore, the proof writer can reduce the cognitive cost of making connections by paying close attention to how the proof \emph{manages information}.\footnote{\citet{AvigadModularity} has discussed information management in the form of modularity and \citet{Sieg2013} has discussed the related concept of ``hierarchical organization.''  Additionally \citet{Author2014,Author2016} have analyzed a detailed case study from the history of mathematics highlighting different forms of information management and their importance.}  Generally speaking, a proof that manages information well highlights information when it is relevant and hides it when it is not.  The proof writer should consider how the proof manages information at both a local and global level.  At the local level, for example, displaying too much information can make the proof much harder to parse than it needs to be.  Making the proof harder to parse means that it will take more effort for the reader to make connections between her context and the proof, and thus make it less likely that the steps will be role and insight motivated.

As a concrete example, consider Wilson's Theorem.\footnote{For a detailed historical case study on Wilson's Theorem and motivation, see \citep[\S 3.4]{Author2015}.}  This result from elementary number theory states that if $p$ is prime, then $(p-1)! + 1$ is evenly divisible by $p$.  The first proof of this theorem was given by Lagrange and works by considering the polynomial $(x+1)(x+2)\ldots(x+p-1)$ and setting up a recurrence relation between its coefficients.  The recurrence relation is exploited to show that the constant coefficient plus one is evenly divisible by $p$ and the constant coefficient is then shown to be $(p-1)!$.  However, different ways of presenting the proof can have a big effect on how easy it is to parse.

Consider, for example, the proof that starts as follows:

\begin{quote}
\textbf{Proof Extract (Style 1)}
\[\mbox{Let } (x+1)(x+2)\ldots(x+p-1) = x^{p-1}+A_{1}x^{p-2}+\ldots+A_{p-1}\]
   Then we have
     \[(x+2)(x+3)\ldots(x+p) = (x+1)^{p-1}+A_{1}(x+1)^{p-2} +\ldots+A_{p-1}\]
   And so we can infer
   
    \[(x+p)(x^{p-1}+A_{1}x^{p-2}+\ldots+A_{p-1})\]
    \[= (x+1)^{p} + A_{1}(x+1)^{p-1}+\ldots+A_{p-1}(x+1)\]
        \ldots
\end{quote}

Compare this with the following, which presents the same reasoning in a different way:

\begin{quote}
\textbf{Proof Extract (Style 2)}
 \[\mbox{Let } L(x) = (x+1)(x+2)\ldots(x+p-1)\]
    Then we have
    \[L(x+1) = (x+2)(x+3)\ldots(x+p)\]
    And so we can infer\[
        (x+p)L(x)=(x+1)L(x+1)\]
        \ldots
\end{quote}

The second style of proof hides information when it is not relevant.  More precisely, it does not introduce the coefficients $A_{1}, A_{2}, \ldots, A_{p-1}$ at the start of the proof, as they are not yet needed.  Removing this information, and including a handy abbreviation for the polynomial, makes the resulting proof extract much easier to parse.

However, how the proof manages information at a more global level is also important.  In particular, if the structure of a proof is clear, then this puts the reader in a good position to make relevant connections between the proof and her context.  For example, if the proof is split up into a series of lemmas, then seeing that step $S$ occurs in the proof of lemma $L$ can help direct the reader's search by narrowing down the connections she should consider.  This in turn can help her to more easily identify the intended task(s) of $S$, as well as where it reasonably came from.

Breaking out lemmas\footnote{Breaking out lemmas is one way to increase the modularity of a proof.  \citet{AvigadModularity} has analyzed the benefits of such proofs.  Note that if we break out lemmas, then we must decide when and how to present these lemmas to the readers.  These issues relate directly to Sieg's notion of ``hierarchical organization'' \citep{Sieg2013}.}  from within a large complicated proof is thus one way that a proof writer can help manage information effectively on a more global scale.  Another possible approach for complicated proofs is to provide a brief sketch or outline.  This will again provide the reader with information about the proof structure that can help guide her search for connections and thus help her to more easily identify each step's intended task(s), as well as where they reasonably came from.

Ultimately, writing motivated proofs relative to a given context can be difficult.  Indeed, historical case studies document the substantial time and effort required to craft proofs which manage information efficiently.\footnote{See, for example, the case study discussed by \citep{Author2014,Author2016} in which such changes took around 100 years.}  Nonetheless, because motivated proofs help disseminate and advance mathematical knowledge, it is time and effort well spent!

\section{Motivation and Other Virtues}\label{othervirtues}
Now that we have obtained a more precise analysis of motivated proofs, we can consider the relationship between motivational power and other desirable proof qualities.  Here I will focus on explanatory power, beauty and fit.

\subsection{Explanation}
The explication of motivated proofs, proposed in section (\ref{Def}), may invite comparisons to explanatory proofs.  For example, we might think that role and insight motivation are types of explanation, and so conclude that a motivated proof is also explanatory.  However, it is difficult to compare motivated proofs to explanatory ones because there is little consensus over what constitutes an explanatory proof.  Indeed, mathematicians and philosophers have strong and vastly different intuitions about which proofs are explanatory.\footnote{See e.g. \citep{Lange} for a discussion of competing intuitions over proofs by induction. }  Moreover, the accounts of explanatory proofs that have been proposed are also strikingly different from each other.  For example, \citet{steiner1978} focused on properties of individual proofs, while \citet{Kitcher1989} focused on unification of a mathematical domain.

Nonetheless, while it is difficult to make a full comparison, there do seem to be important differences between motivated proofs and explanatory ones.  Elsewhere \citep[\S 5]{Author2015} I have proposed certain proofs that are, according to my account, motivated with respect to specific natural contexts but fail to be explanatory according to Steiner's and Kitcher's theories.  At a more general level, we can also note that theories of mathematical explanation are not usually sensitive to whether the explanation is \emph{recognized} as such by mathematicians.  However, such recognition is built into the definition of role and insight motivation---recall that for a step to count as role or insight motivated relative to a given context, agents with access to that context must \emph{identify} the intended task of each step, as well as where each step could have reasonably come from.

Moreover there seem to be important differences between the instrumental value of motivated and explanatory proofs. Elsewhere \citep{Author2019} I identified reuse of resources as an important mathematical goal and argued that explanatory proofs do not generally help mathematicians to achieve it. However, as we have seen in section (\ref{value}), motivated proofs communicate how and when to apply mathematical resources and thus do promote reuse of their resources.

Nonetheless, there may be a connection between motivated proofs and \emph{scientific} explanation.\footnote{I am grateful to an anonymous reviewer for pointing out this connection.}  According to Hempel's deductive-nomological account, an explanation shows that the occurrence of some phenomenon was to be expected, given the specific conditions and relevant laws \citep{Woodward2017-ql}.  The notion of insight motivation involves identifying where a given step could have reasonably come from and thus seems to relate to grasping why that step was to be expected.

The relationship between motivational and explanatory power is thus both interesting and complicated and should be further explored.

\subsection{Beauty}
Although the notion of beauty, like explanation, provokes conflicting intuitions, there is a closer connection between motivational power and beauty, at least on Cellucci's account.  According to Cellucci ``\ldots a mathematical demonstration or theorem is beautiful when it provides understanding'' \citep[abstract]{Cellucci}, where ``understanding'' is defined as ``recognition of the fitness of the parts to each other and to the whole'' \citep[\S 8, para 3]{Cellucci}.  Cellucci further explained: ``There is fitness of the parts of a demonstration to each other and to the whole when it is clear what the whole idea of the demonstration is, what the contribution of each part of the demonstration to the whole idea is, and why such contribution is essential'' \citep[\S 10, para 1]{Cellucci}.

The notions of role and insight motivation seem to capture some of the ways in which proof steps can ``fit'' with each other and the whole demonstration in the way that Cellucci described.  Further, as readers of a motivated proof, relative to a given context, can identify the intended task(s) of each step, as well as where each step could have reasonably come from, they would thus recognize their fit.  Consequently, it appears that a motivated proof, relative to a suitable natural context, is also beautiful on Cellucci's account.

\subsection{Fit}
We have already seen, in the discussion of beauty, that there is a close relationship between motivated proofs and ``fit.''  \citet{Fit2016} have analyzed fit in its own right and have identified three general kinds, each of which has two subtypes: direct fit (coherence, specificity); presentational fit (level of detail, transparency); and familial fit (generality, connectedness).  Their notion of transparency seems most related to motivated proofs.  Indeed, Raman-Sundstr\"{o}m and \"{O}hman described transparency in the following terms ``In a proof that is strong in this criterion, it is easy to see `what is going on.'  In other words, the structure of the proof is natural for the particular argument and there is no \emph{deus ex machina} component'' \citep[5]{Fit2016}.  They further explained that ``\ldots if a proof is transparent, a reader with the appropriate background should be in an ideal position to grasp the ideas of the proof'' \citep[5]{Fit2016}.

As readers of a motivated proof, relative to a suitable context, can identify each step's intended task(s), as well as where each step could have reasonably come from, it should indeed be easy for them to grasp ``what is going on.''  Consequently, motivated proofs also appear to be transparent.  However, perhaps it is possible to grasp ``what is going on'' in a proof without all of its steps being role or insight motivated.  For example, perhaps it is possible to understand the core ideas of a proof and have a high level grasp of ``what is going on'' within it while still finding a few steps involving the technical details to be puzzling.  If so, then it is possible for there to be proofs which are transparent but fail to be fully motivated.

\section{Concluding Remarks}\label{conclude}
In this paper I have offered a context-sensitive explication of motivated proofs.  Recall from section (\ref{Def}) that a proof is motivated relative to a given context if and only if each step is both role and insight motivated relative to the context.  In other words, a proof is motivated relative to a given context if and only if a reader with access to that context can identify the tasks each step was intended to perform, as well as where each step could have reasonably come from.  This explication thus reflects the intuition that motivated proofs are those which do not contain any ``puzzling'' steps, as role and insight motivation rule out potential sources of confusion for the reader.  Moreover, as we saw in section (\ref{value}), such proofs have three main epistemic benefits: (i) they promote understanding; (ii) they successfully convey new mathematical resources; (iii) they stimulate new discoveries.  Ultimately, then, such proofs serve to disseminate and advance mathematical knowledge among practitioners of a given context.

However, there may be other, alternative definitions of motivated proofs that are of interest. For example, while the account I have developed here is reductionist, it may be possible to give an alternative, holistic, definition focusing on how interesting a proof is or how it connects to other areas of mathematics.  Further, mathematicians often speak of proofs which are, for example, ``geometrically motivated'' \citep{Sherali} or ``physically motivated'' \citep{Halliwell}.  Proofs which are motivated in these ways are worthy of investigation in their own right.

Moreover, while I have focused on motivated \emph{proofs} here, other mathematical artifacts are often said to be motivated.  For example, mathematicians call attention to definitions and theories which are (or fail to be) well motivated.  Here are two examples:

\begin{quote}
It is well known that not all algorithms are feasible; whether an algorithm is feasible or not depends on how many computational steps this algorithm requires.  The problem with the existing definitions of feasibility is that they are rather \emph{ad hoc}.  Our goal is to use the maximum entropy (MaxEnt) approach and get more motivated definitions \citep[25]{Cooke98}
\end{quote}

\begin{quote}
Starting from a small number of well-motivated axioms, we derive a unique definition of sums with a noninteger number of addends. \citep{Muller11}
\end{quote}

This suggests that it may be of interest to develop an account of motivational power that applies to these artifacts as well.  Indeed, there will likely be connections between the accounts themselves, as well as between accounts of other virtues.  I have already discussed, in section (\ref{othervirtues}), the relationship between motivated proofs, explanatory power, beauty and fit.  Developing accounts of motivated definitions and theories may reveal further connections to virtues such as fruitfulness, which has been analyzed by \citet{FandN} and \citet{Yap2011}.  Finally, given the close connection between mathematics and science, the notions of motivated mathematical artifacts may generalize so as to apply to scientific concepts and theories.  Consequently the work presented here is just the beginning of an investigation into motivational power.

However, it is not just philosophical work that remains to be done.  Investigation into the topic of motivated proofs, and motivated mathematical and scientific artifacts more generally, should be highly interdisciplinary, incorporating insights from, for example, philosophy, history, psychology and education.

\vspace{0.5cm}
\noindent \textbf{Acknowledgments}
I would like to thank Jeremy Avigad, Michael Detlefsen, Michael Friedman, Yacin Hamami, Erich Kummerfeld, Kenneth Manders, Gihan Marasingha, Wilfried Sieg and the anonymous reviewers for their helpful discussions and feedback. I am also grateful to audiences at conferences and seminars in Fullerton, USA, Nancy, France and Brussels, Belgium for their helpful questions and comments.  Part of this work was undertaken while I held a Postdoctoral Scholarship at the Suppes Center for History and Philosophy of Science at Stanford University.

\bibliographystyle{plainnat}
\bibliography{Motivation}

\begin{thebibliography}{37}
\providecommand{\natexlab}[1]{#1}
\providecommand{\url}[1]{\texttt{#1}}
\expandafter\ifx\csname urlstyle\endcsname\relax
  \providecommand{\doi}[1]{doi: #1}\else
  \providecommand{\doi}{doi: \begingroup \urlstyle{rm}\Url}\fi

\bibitem[{Art of Problem Solving}(n.d.)]{AoPS}
{Art of Problem Solving}.
\newblock Proofs of {A}{M}-{G}{M}, n.d.
\newblock URL
  \url{https://artofproblemsolving.com/wiki/index.php/Proofs_of_AM-GM}.
\newblock Archived at \url{https://perma.cc/QTE7-YEH8}.

\bibitem[Avigad(2018)]{AvigadModularity}
Jeremy Avigad.
\newblock {Modularity} {in} {mathematics}.
\newblock \emph{Review of Symbolic Logic}, pages 1--33, February 2018.

\bibitem[Avigad and Morris(2014)]{Author2014}
Jeremy Avigad and Rebecca~Lea Morris.
\newblock The concept of ``character'' in {D}irichlet's theorem on primes in an
  arithmetic progression.
\newblock \emph{Archive for history of exact sciences}, 2014.

\bibitem[Avigad and Morris(2016)]{Author2016}
Jeremy Avigad and Rebecca~Lea Morris.
\newblock Character and object.
\newblock \emph{Review of Symbolic Logic}, 2016.

\bibitem[Cellucci(2015)]{Cellucci}
Carlo Cellucci.
\newblock Mathematical beauty, understanding, and discovery.
\newblock \emph{Foundations of science}, 20\penalty0 (4):\penalty0 339--355,
  November 2015.

\bibitem[Clark(2013)]{Clark2Squares}
Peter~L. Clark.
\newblock Proof that every odd integer is a difference of two squares.
\newblock URL \url{http://math.stackexchange.com/q/510239}. Archived at
  \url{https://perma.cc/F8UZ-G7SP}., 2013.

\bibitem[Cohn(2006)]{CFE}
Henry Cohn.
\newblock A short proof of the simple continued fraction expansion of $e$.
\newblock \emph{The American Mathematical Monthly}, 113\penalty0 (1):\penalty0
  57--62, 2006.

\bibitem[Cooke et~al.(1998)Cooke, Kreinovich, and Longpr{\'e}]{Cooke98}
D~E Cooke, V~Kreinovich, and L~Longpr{\'e}.
\newblock Which algorithms are feasible? {M}axent approach.
\newblock In Gary~J Erickson, Joshua~T Rychert, and C~Ray Smith, editors,
  \emph{Maximum Entropy and Bayesian Methods}, Fundamental Theories of Physics,
  pages 25--33. Springer Netherlands, 1998.

\bibitem[Corr{\'a}di and Szad{\'o}(1993)]{HT}
Kereszt{\'e}ly Corr{\'a}di and S{\'a}ndor Szad{\'o}.
\newblock A generalized form of {H}aj{\'o}s' theorem.
\newblock \emph{Communications in Algebra}, 21\penalty0 (11):\penalty0
  4119--4125, January 1993.

\bibitem[Detlefsen and Arana(2011)]{Purity}
Michael Detlefsen and Andrew Arana.
\newblock Purity of methods.
\newblock \emph{Philosophers' Imprint}, 11, 2011.

\bibitem[Deutsch(2012)]{deutsch}
Frank~R Deutsch.
\newblock \emph{Best Approximation in Inner Product Spaces}.
\newblock Springer Science \& Business Media, December 2012.

\bibitem[Gowers(2008)]{ReduceToSimple}
Timothy Gowers.
\newblock Very brief tricki update, 2008.
\newblock URL
  \url{https://gowers.wordpress.com/2008/11/28/very-brief-tricki-update/}.
\newblock Archived at \url{https://perma.cc/9V4P-734F}.

\bibitem[Halliwell(2014)]{Halliwell}
J~J Halliwell.
\newblock Two proofs of fine's theorem.
\newblock \emph{Physics letters. A}, 378\penalty0 (40):\penalty0 2945--2950,
  August 2014.

\bibitem[Jones(2010)]{IPI}
Timothy~W Jones.
\newblock Discovering and proving that $\pi$ is irrational.
\newblock \emph{The American Mathematical Monthly}, 117\penalty0 (6):\penalty0
  553--557, June 2010.

\bibitem[Kitcher(1989)]{Kitcher1989}
Philip Kitcher.
\newblock Explanatory unification and the causal structure of the world.
\newblock In Philip Kitcher and Wesley Salmon, editors, \emph{Scientific
  Explanation}, pages 410--505. Minneapolis: University of Minnesota Press,
  1989.

\bibitem[Lange(2009)]{Lange}
Marc Lange.
\newblock Why proofs by mathematical induction are generally not explanatory.
\newblock \emph{Analysis}, 69\penalty0 (2):\penalty0 203--211, April 2009.

\bibitem[Lange(2015)]{LangeDepth2015}
Marc Lange.
\newblock Depth and explanation in mathematics.
\newblock \emph{Philosophia Mathematica. Series III}, 23\penalty0 (2):\penalty0
  196--214, June 2015.

\bibitem[MacLane(1935)]{MacLane1935}
Saunders MacLane.
\newblock A logical analysis of mathematical structure.
\newblock \emph{The Monist}, 45\penalty0 (1):\penalty0 118--130, 1935.

\bibitem[Montanaro(2013)]{InequalityProofs}
Ashley Montanaro.
\newblock Simplest or nicest proof that $1+x \le e^x$.
\newblock Mathematics Stack Exchange, 2013.
\newblock URL \url{https://math.stackexchange.com/q/504663}.
\newblock Archived at \url{https://perma.cc/E3Y5-WW9F}.

\bibitem[Morris(2015)]{Author2015}
Rebecca~Lea Morris.
\newblock \emph{Appropriate Steps: A Theory of Motivated Proofs}.
\newblock PhD thesis, Carnegie Mellon University, 2015.

\bibitem[Morris(2019)]{Author2019}
Rebecca~Lea Morris.
\newblock Do mathematical explanations have instrumental value?
\newblock \emph{Synthese}, February 2019.

\bibitem[M{\"u}ller and Schleicher(2011)]{Muller11}
Markus M{\"u}ller and Dierk Schleicher.
\newblock How to add a noninteger number of terms: from axioms to new
  identities.
\newblock \emph{The American Mathematical Monthly}, 118\penalty0 (2):\penalty0
  136--152, 2011.

\bibitem[P\'{o}lya(1949)]{WWM}
George P\'{o}lya.
\newblock With, or without, motivation?
\newblock \emph{The American Mathematical Monthly}, 56\penalty0 (10):\penalty0
  684--691, 1949.

\bibitem[Raman-Sundstr{\"o}m and {\"O}hman(2018)]{Fit2016}
Manya Raman-Sundstr{\"o}m and Lars-Daniel {\"O}hman.
\newblock Mathematical fit: A case study.
\newblock \emph{Philosophia Mathematica}, 26\penalty0 (2):\penalty0 184--210,
  June 2018.

\bibitem[Rav(1999)]{Rav1999}
Yehuda Rav.
\newblock Why do we prove theorems?
\newblock \emph{Philosophia Mathematica}, 7\penalty0 (1):\penalty0 5--41,
  February 1999.

\bibitem[Rogers(1973)]{rogers_lusins_1973}
C~A Rogers.
\newblock Lusin's second separation theorem.
\newblock \emph{Journal of the London Mathematical Society. Second Series},
  s2-6\penalty0 (3):\penalty0 491--503, May 1973.

\bibitem[Rota(1997)]{Rota}
Gian-Carlo Rota.
\newblock The phenomenology of mathematical beauty.
\newblock \emph{Synthese}, 111\penalty0 (2):\penalty0 171--182, May 1997.

\bibitem[Sandborg(1997)]{SandborgDiss}
David Sandborg.
\newblock \emph{Explanation in Mathematical Practice}.
\newblock PhD thesis, University of Pittsburgh, 1997.

\bibitem[Sandborg(1998)]{Sandborg1998}
David Sandborg.
\newblock Mathematical explanation and the theory of why-questions.
\newblock \emph{The British journal for the philosophy of science}, 49\penalty0
  (4):\penalty0 603--624, December 1998.

\bibitem[Sherali(1987)]{Sherali}
Hanif~D Sherali.
\newblock A constructive proof of the representation theorem for polyhedral
  sets based on fundamental definitions.
\newblock \emph{American Journal of Mathematical and Management Sciences},
  7\penalty0 (3-4):\penalty0 253--270, February 1987.

\bibitem[Sieg(2010)]{Sieg2013}
Wilfried Sieg.
\newblock Searching for proofs (and uncovering capacities of the mathematical
  mind).
\newblock \emph{Proofs, Categories and Computations--Essays in Honor of Grigori
  Mints}, pages 189--215, 2010.

\bibitem[Steele(2004)]{CSMasterclass}
J~Michael Steele.
\newblock \emph{The {Cauchy-Schwarz} Master Class : An Introduction to the Art
  of Mathematical Inequalities}.
\newblock Cambridge University Press, 2004.

\bibitem[Steiner(1978)]{steiner1978}
Mark Steiner.
\newblock Mathematical explanation.
\newblock \emph{Philosophical studies}, 34\penalty0 (2):\penalty0 135--151,
  August 1978.

\bibitem[Tappenden(2008)]{FandN}
Jamie Tappenden.
\newblock Mathematical concepts: Fruitfulness and naturalness.
\newblock In Paolo Mancosu, editor, \emph{The Philosophy of Mathematical
  Practice}, pages 276--301. Oxford University Press, 2008.

\bibitem[Thurston(1994)]{Thurston}
William~P Thurston.
\newblock On proof and progress in mathematics.
\newblock \emph{Bulletin of the American Mathematical Society}, 30:\penalty0
  161--177, 1994.

\bibitem[Woodward(2017)]{Woodward2017-ql}
James Woodward.
\newblock Scientific explanation.
\newblock In Edward~N Zalta, editor, \emph{The Stanford Encyclopedia of
  Philosophy}. Metaphysics Research Lab, Stanford University, fall 2017
  edition, 2017.
\newblock URL
  \url{https://plato.stanford.edu/archives/fall2017/entries/scientific-explanation/}.

\bibitem[Yap(2011)]{Yap2011}
Audrey Yap.
\newblock Gauss' quadratic reciprocity theorem and mathematical fruitfulness.
\newblock \emph{Studies in History and Philosophy of Science. Part B. Studies
  in History and Philosophy of Modern Physics}, 42\penalty0 (3):\penalty0
  410--415, September 2011.

\end{thebibliography}

\end{document}